\documentstyle[12pt]{article}
\def\Bbb R{{\rm \bf R}}
\def\proclaim#1{\vskip2mm{\bf #1}\em}
\def\endproclaim{\em \vskip2mm}
\def\tag#1{\eqno(#1)}
\def\gathered{\begin{array}{c}}
\def\endgathered{\end{array}}
\def\text{\mbox}

\begin{document}

\title {On reconstruction of inclusions in a heat conductive body from dynamical
boundary data over a finite time interval}
\author{Masaru IKEHATA\footnote{
Department of Mathematics,
Graduate School of Engineering,
Gunma University, Kiryu 376-8515, JAPAN}
and
Mishio KAWASHITA\footnote{
Department of Mathematics,
Graduate School of Sciences,
Hiroshima University, Higashi Hiroshima 739-8526, JAPAN}}
\date{July 1, 2010 Final}
\maketitle
\begin{abstract}
The {\it enclosure method} was originally introduced for inverse problems
of concerning non-destructive evaluation governed by elliptic equations.
It was developed as one of useful approaches in inverse problems and applied for various equations.
In this paper, an application of the enclosure method to an inverse initial boundary value problem
for a parabolic equation with a discontinuous coefficient is given.
A simple method to extract the {\it depth} of
unknown inclusions in a heat conductive body from a {\it single set} of
the temperature and heat flux on the boundary observed over a
{\it finite} time interval is introduced.  Other related results with {\it infinitely
many} data are also reported.  One of them gives the {\it minimum radius}
of the open ball centered at a given point that contains the inclusions.
The formula for the minimum radius is newly discovered.

\noindent
AMS: 35R30, 80A23

\end{abstract}


\section{Introduction}
A part in the body that has a different conductivity from the
known reference one is called an {\it inclusion}.
In this paper, we consider a mathematical formulation of
an inverse problem for the heat equation with inclusion.
Assume that we have a set of the pair of the temperature field on
the boundary of a heat conductive body and the corresponding heat
flux across the boundary of the body over a finite time interval.
This is the set of observed data in our inverse problem.
The problem in this paper is to find
what information about inclusions in the body can one
extract from the set.
The solution to this problem may have possible application
to non-destructive evaluation by
{\it thermal imaging}.
We study this problem from a mathematical
point of view and aim at seeking an analytical approach for extracting
information about the location and shape of the inclusions.

Let $\Omega$ be a bounded domain of $\Bbb R^3$ with a smooth boundary.
We denote the unit outward normal vectors to $\partial\Omega$ by the symbol $\nu$.
Let $T$ be an arbitrary {\it fixed} positive number.

Given $f=f(x,t),\,(x,t)\in\partial\Omega\times\,]0,\,T[$
let $u=u(x,t)$ be the solution of the initial boundary
value problem for the parabolic equation
$$
\begin{array}{c}
\displaystyle u_t-\nabla\cdot\gamma\nabla u=0\,\,\text{in}\,\Omega\times\,]0,\,T[,\\
\\
\displaystyle
\gamma\nabla u\cdot\nu=f\,\,\text{on}\,\partial \Omega\times\,]0,\,T[,\\
\\
\displaystyle
u(x,0)=0\,\,\text{in}\,\Omega,
\end{array}
\tag {1.1}
$$
where $\gamma=\gamma(x)=(\gamma_{ij}(x))$ satisfies

$\quad$

\noindent
(G1)  for each $i,j=1,2,3$ $\gamma_{ij}(x)=\gamma_{ji}(x)\in L^{\infty}(\Omega)$;

$\quad$

\noindent
(G2)  there exists a positive constant $C$ such that $\gamma(x)\xi\cdot\xi\ge C\vert\xi\vert^2$ for all $\xi\in\Bbb R^3$
and a. e. $x\in\Omega$.

$\quad$

This paper is concerned with the extraction of information about
`discontinuity' of $\gamma$ from $u$ and $\gamma\nabla
u\cdot\nu$ on $\partial\Omega\times ]0,\,T[$ for some $f$ and an
arbitrary fixed $T<\infty$. However, we do not consider completely
general $\gamma$.
Instead we assume that there exists an open set
$D$ with a smooth boundary such that $\overline D\subset\Omega$
and $\gamma(x)$ a.e. $x\in\Omega\setminus D$ coincides with
the $3\times 3$ identity matrix $I_3$ and satisfies one of the
following two conditions:

$\quad$

\noindent
(A1)  there exists a positive constant $C'$ such that $-(\gamma(x)-I_3)\xi\cdot\xi\ge C'\vert\xi\vert^2$
for all $\xi\in\Bbb R^3$ and a.e. $x\in D$;

$\quad$

\noindent
(A2)  there exists a positive constant $C'$ such that $(\gamma(x)-I_3)\xi\cdot\xi\ge C'\vert\xi\vert^2$
for all $\xi\in\Bbb R^3$ and a.e. $x\in D$.

$\quad$

Write $h(x)=\gamma(x)-I_3$ a.e. $x\in D$.
In this paper we consider the following problem:

$\quad$

{\bf\noindent Inverse Problem.}  Assume that both $D$ and $h$ are {\it unknown}.
Extract information about the location and shape
of $D$ from a set of the pair of temperature $u(x,t)$ and
heat flux $f(x,t)$ for $(x,t)\in\,\partial\Omega\,\times\,]0,\,T[$.

$\quad$

$D$ is a model of the union of unknown {\it inclusions} where
the heat conductivity is anisotropic, different from that of the surrounding homogeneous isotropic conductive medium.
The problem is a mathematical formulation of a typical inverse problem in thermal
imaging.

Elayyan-Isakov \cite{EI} investigated the {\it uniqueness issue}
of this type of problem. As a corollary of their uniqueness
theorem we know that the lateral Neumann-to-Dirichlet map:
$f\longmapsto u\vert_{\partial\Omega\times]0,\,T[}$ uniquely
determines $D$ together with $h$ inside $D$ if
$\Omega\setminus\overline D$ is connected and $h$ is given by $b
I_3$ with a smooth function $b$ on $\overline D$. However, their
purpose is to recover the {\it full information} about the
location and shape of $D$ and for the purpose their proof requires
{\it infinitely many pairs} of the temperature and heat flux on
$\partial\Omega\times ]0,\,T[$ even just for determining a single
point on $\partial D$. This shows a difficulty of obtaining the
detailed image of inclusions from boundary measurements.

Also note that in \cite{DKN} an approach to Inverse Problem
in a one-space dimensional case which is an analogue of the {\it
probe method} introduced by Ikehata \cite{IP1} has been proposed.
Therein $\Omega$ and $D$ are given by open intervals. Their
approach introduces the {\it pre-indicator function} denoted by
$I(y',s';y,s)$ with their notation which depends on {\it two
points}  $(y,s), (y',s')\in (\Omega\setminus\overline D)\times\,]0,\,T[$.
The pre-indicator function can
be computed from an integral over the lateral boundary
$\partial\Omega\times\,]0,\,T[$ which involves the
Neumann-to-Dirichlet map acting on {\it infinitely many} special
heat fluxes prescribed on the lateral boundary.
The procedure in \cite{DKN} for the determination of
the right-end point of $D$ consists of two steps: (i) computation
of $I(y+\epsilon,s+\epsilon^2;y,s)$ for each $\epsilon>0$ at given point $(y,s)$ in the
spacetime and
$\lim_{\epsilon\longrightarrow
0}I(y+\epsilon,s+\epsilon^2;y,s)$; (ii) moving $(y,s)$ in (i) to left along a path on the
spacetime and do the same computation until the computation
result becomes large.
This procedure has been tested numerically in \cite{DLLN}.

In this paper we mainly seek a simpler method that yields a {\it partial} or {\it rough} information about the location and shape
of $D$ from $u(x,t)$ on $\partial\Omega\times\,]0,\,T[$ for a single {\it fixed} heat flux or
{\it explicit} heat fluxes prescribed on $\partial\Omega\times\,]0,\,T[$.
We think that
this type of study provides us knowledge about {\it good} heat fluxes on the boundary of the body
to get such information.
In \cite{IK2} we have already developed an argument based on the
{\it enclosure method} which was originally introduced for elliptic
equations in \cite{IS,IE}
to derive two types of formulae in the case
when the inclusion has the {\it zero conductivity}, that is a {\it
cavity}. The argument yields the values of the support function of
the cavity at a given direction and the distance of a given point
outside the body to the cavity from the temperature fields and
special explicit heat fluxes.
In this paper, we see that the argument also
works for the inclusion case and also yields new information:
the {\it minimum radius} of the open ball centered at a given point that contains
the inclusions.

The main new point of this paper is as follows: an introduction of another argument which is also based on the enclosure method
and gives a formula which has not been considered in \cite{IK2}.
It makes use of a single set of a general heat flux and the corresponding
temperature field on the surface of the body over a {\it finite time interval}.
It yields a {\it depth} of unknown inclusions in a heat conductive body from the surface of the body.
We give a heat flux $f(x, t)$ on the surface of the body
as it moves from the beginning of the experiments (cf (1.2) below).
But, there is no other condition on $f$.
Hence note that we do not need to prescribe any explicit heat flux.

Note that in Theorem 2.1 of \cite{I4} the enclosure method has been applied
to a one-space dimensional version of Inverse Problem.  Therein
{\it complex} exponential solutions of the backward heat equation with a large parameter are used.
In this paper we use {\it only} real exponential solutions.

\subsection{A formula with a general heat flux}

The new point of this paper is a derivation of the following
formula which can be considered
as the main result of this paper.
It makes use of a single set of a heat flux and the corresponding temperature
on $\partial\Omega\times\,]0,\,T[$ and gives partial or rough information about the location
of inclusions.

\proclaim{\noindent Theorem 1.1.} Let $f=f(x,t)$, being square integrable in $(x,t)$, satisfy: there
exists $\mu\in\Bbb R$ such that
$$\displaystyle
0<\inf_{x\in\partial\Omega}\liminf_{\tau\longrightarrow\infty}\tau^{\mu}\int_0^Te^{-\tau t}f(x,t)dt\le
\sup_{x\in\partial\Omega}\limsup_{\tau\longrightarrow\infty}\tau^{\mu}
\int_0^Te^{-\tau t}f(x,t)dt<\infty
\tag {1.2}
$$
and that the function
$$\displaystyle
\partial\Omega\ni x\longmapsto g(x;\tau)\equiv\int_0^Te^{-\tau t}f(x,t)dt
$$
is continuous.

Let $u=u_f(x,t)$ be the weak solution of (1.1) for this $f$
and let $v=v_g(x;\tau)$ be the solution of
$$\begin{array}{c}
\displaystyle
(\triangle-\tau)v=0\,\,\text{in}\,\Omega,\\
\\
\displaystyle
\frac{\partial v}{\partial\nu}
=g(x;\tau)\,\,\text{on}\,\partial\Omega.
\end{array}
\tag {1.3}
$$
Then, there exists $\tau_0>0$ such that

$\bullet$  if (A1) is satisfied, then for all $\tau\ge\tau_0$
$$\displaystyle
\int_{\partial\Omega}\int_0^Te^{-\tau t}\left(v_g(x;\tau)f(x,t)-u_f(x,t)g(x;\tau)\right)dtdS<0;
$$

$\bullet$  if (A2) is satisfied, then for all $\tau\ge\tau_0$
$$\displaystyle
\int_{\partial\Omega}\int_0^Te^{-\tau t}\left(v_g(x;\tau)f(x,t)-u_f(x,t)g(x;\tau)\right)dtdS>0.
$$
In both cases the formula
$$\displaystyle
\lim_{\tau\longrightarrow\infty}
\frac{1}{2\sqrt{\tau}}
\log\left\vert\int_{\partial\Omega}
\int_0^Te^{-\tau t}\left(v_g(x;\tau)f(x,t)-u_f(x,t)g(x;\tau)\right)dtdS\right\vert
=-\text{dist}\,(D,\partial\Omega)
\tag {1.4}
$$
is valid, where
$$\displaystyle
\text{dist}\,(D,\partial\Omega)=\inf\{\vert y-x\vert\,\vert\,y\in\partial\Omega,\,x\in\,D\}.
$$

\endproclaim

Simple examples of $f(x,t)$ are $f(x,t)=\varphi(t)h(x)$ where $h\in C(\partial\Omega)$ with a positive lower bound
and $\varphi(t)=t^n, n=0,1,2,3,\cdots$ for $0<t<T$.
This is because of
$$
\displaystyle
\lim_{\tau\longrightarrow\infty}\tau^{n+1}\int_0^{T}e^{-\tau t}t^ndt
=\int_0^{\infty}e^{-\xi}\xi^nd\xi>0
$$
and thus (1.2) is satisfied for $\mu=n+1$.
It is also possible to give more general example.

Note that Varadhan \cite{V} considered the asymptotic behaviour as $\tau\longrightarrow\infty$ of
the solution of the problem
$$\begin{array}{c}
\displaystyle
(\triangle-\tau)v=0\,\,\text{in}\,\Omega,\\
\\
\displaystyle
v
=1\,\,\text{on}\,\partial\Omega.
\end{array}
$$
He used the behaviour to establish the short-time asymtotics of
the {\it heat kernel}.  See also \cite{N} and references therein for the subject itself.
Theorem 1.1 shows that this type of solutions can
be applied to inverse initial boundary value problems for parabolic equations
over a {\it finite} time interval.

In \cite{I} Ikehata considered an inverse obstacle scattering problem whose governing equation
is given by the wave equation in three dimensions.
The observation data are given by a wave
field measured on a known surface surrounding unknown obstacles
over a {\it finite} time interval. The wave is generated by an
initial data with compact support outside the surface.
Applying the idea of the enclosure method,
he established an extraction formula of the distance from a given point outside
the surface to obstacles from the data.
To establish the formula he made use of the solution $v\in H^1(\Bbb R^3)$ of
the inhomogeneous modified Helmholtz equation
$$\displaystyle
(\triangle-\tau^2)v+f(x)=0\,\,\text{in}\,\Bbb R^3,
$$
where $f(x)$ is an initial data of the wave field.  Thus the
equations in (1.3) correspond to this equation.  However, in
contrast to the solution of this equation, that of (1.3) has not an
explicit form in general.  In this paper, we solve (1.3) by using
the potential theory and study its behaviour as
$\tau\longrightarrow\infty$ to get a necessary estimate.

\subsection{Other three formulae with special heat fluxes}

If one uses special heat fluxes, then one can explicitly obtain
more information about the location and shape of $D$.  The idea
for the derivation of the following formulae come from \cite{IK2}.

The second result is the following.
\proclaim{\noindent Theorem 1.2.}
Given $\omega\in S^{2}$ let $f$ be the function of $(x,t)\in\partial\Omega\times]0,\,T[$ having a parameter $\tau>0$
defined by the equation
$$\displaystyle
f(x,t;\tau)=\frac{\partial v}{\partial\nu}(x;\tau)\varphi(t),
\tag {1.5}
$$
where $v(x;\tau)=e^{\sqrt{\tau}x\cdot\omega}$ and $\varphi\in L^2(0,\,T)$
satisfying the following condition:
there exists $\mu\in\Bbb R$ such
that
$$\displaystyle
\liminf_{\tau\longrightarrow\infty}
\tau^{\mu}
\left\vert\int_0^Te^{-\tau t}\varphi(t)dt\right\vert>0.
\tag {1.6}
$$
Let $u_f=u_f(x,t)$ be the weak solution of (1.1) for $f=f(x,t;\tau)$
and $h_D(\omega)=\sup_{x\in D}x\cdot\omega$.
Then the formula
$$\displaystyle
\lim_{\tau\longrightarrow\infty}
\frac{1}{2\sqrt{\tau}}
\log\left\vert
\int_{\partial\Omega}
\int_0^Te^{-\tau t}
\left(v(x;\tau)f(x,t;\tau)-u_f(x,t)\frac{\partial v}{\partial\nu}(x;\tau)\right)
dtdS
\right\vert
=h_D(\omega)
\tag {1.7}
$$
is valid.

\endproclaim

\noindent
Note that if $\varphi(t)$ is smooth on $[0,\,T'[$ with $0<T'\le T$ and $t=0$ is not
a zero point with infinite order of $\varphi(t)$, then (1.6) is
satisfied for an appropriate $\mu>0$.

Next we choose a third solution of the equation $(\triangle-\tau)v=0$ in $\Omega$:
given $p\in\Bbb R^{3}\setminus\overline\Omega$
$$\displaystyle
v(x;\tau)=\frac{e^{-\sqrt{\tau}\vert x-p\vert}}{\vert x-p\vert},\,x\in\Omega.
\tag {1.8}
$$
Using this $v$, we obtain the third formula.

\proclaim{\noindent Theorem 1.3.}
Let $p\in\Bbb R^3\setminus\overline\Omega$ and
replace $v$ of $f$ in (1.5) with (1.8).  Let $u_f=u_f(x,t)$ be the weak solution of (1.1) for this $f=f(x,t;\tau,p)$.
Then assuming (1.6), one has the formula
$$\displaystyle
\lim_{\tau\longrightarrow\infty}
\frac{1}{2\sqrt{\tau}}
\log\left\vert\int_{\partial\Omega}
\int_0^Te^{-\tau t}\left(v(x;\tau)f(x,t;\tau,p)-u_f(x,t)\frac{\partial v}{\partial\nu}(x;\tau)\right)dtdS\right\vert
=-d_D(p),
\tag {1.9}
$$
where $d_D(p)$ denotes the distance from $p$ to $D$:
$$\displaystyle
d_D(p)=\inf\{\vert y-p\vert\,\vert\,y\in D\}.
$$

\endproclaim

Finally we introduce another formula which is also
new and not given in \cite{IK2}.
Let $y\in\Bbb R^3$ be an {\it arbitrary} fixed point.
We choose the function $v$ given by
$$\displaystyle
v(x;\tau)=\left\{
\begin{array}{lr}
\displaystyle
\frac{e^{\sqrt{\tau}\vert x-y\vert}-e^{-\sqrt{\tau}\vert x-y\vert}}{\vert x-y\vert}, & \quad\text{if
$x\in\Bbb R^3\setminus\{y\}$,}\\
\\
\displaystyle 2\sqrt{\tau}, & \quad\text{if $x=y$.}
\end{array}
\right.
\tag {1.10}
$$
Note that the $v(x;\tau)$ is smooth as the function of $x$ and satisfies the modified Helmholtz equation
in the whole space.
Hence we can choose the reference point $y\in\Bbb R^3$ without any restriction.
Note that Theorem 1.3 gives $d_D(p)$; however, we have to take $p\in\Bbb R^3\setminus\overline\Omega$.

\proclaim{\noindent Theorem 1.4.}
Let $y\in\Bbb R^3$ and replace $v$ of $f$ in (1.5) with (1.10).
Let $u_f=u_f(x,t)$ be the weak solution of (1.1) for $f=f(x,t;\tau,y)$.
Then assuming (1.6), one has the formula
$$\displaystyle
\lim_{\tau\longrightarrow\infty}
\frac{1}{2\sqrt{\tau}}
\log\left\vert
\int_{\partial\Omega}
\int_0^Te^{-\tau t}
\left(v(x;\tau)f(x,t;\tau,y)-u_f(x,t)\frac{\partial v}{\partial\nu}(x;\tau)\right)
dtdS
\right\vert
=R_D(y),
\tag {1.11}
$$
where $R_D(y)=\sup_{x\in D}\vert x-y\vert$.

\endproclaim

The above theorem makes use of a smooth solution of the modified
Helmholtz equation that grows every point as
$\tau\longrightarrow\infty$. The function $R_D(y)$, $y\in\Omega$
is a newcomer and gives the minimum radius of the ball centered
at $y$ that contains $D$. Moreover we have the estimate of $D$
from above as
$$\displaystyle
D\subset\cap_{y\in\Bbb R^3}\{x\in\Bbb R^3\,\vert\,\vert x-y\vert<R_D(y)\}.
$$

\subsection{Construction of the paper}
A brief outline of this paper is as follows.
Theorem 1.1 is proved in Subsection 2.3 after formulating the notion of the weak solution
of (1.1) together with a related estimate in subsection 2.1.
The proof is based on an integral
identity which is described in subsection 2.2.
Using the identity, we give an asymptotic representation formula of the
integral
$$\displaystyle
\int_{\partial\Omega}\int_0^Te^{-\tau t}\left(v(x;\tau)f(x,t)-u_f(x,t)\frac{\partial v}{\partial\nu}(x;\tau)\right)dtdS
$$
whose leading term is given by using two Neumann-to-Dirichlet maps
for the operators $\triangle-\tau$ and
$\nabla\cdot\gamma\nabla-\tau$ in $\Omega$. Then with a help of a system of
integral inequalities \cite{I} which is widely used in previous
applications of the enclosure method to elliptic equations
\cite{IE} we see that the problem is reduced to giving some
asymptotic estimates for the integral of the gradient of $v_g$
over $D$.  In some sense, this is an {\it indirect} verification
of the hypothesis: $v_g(x;\tau)\sim
e^{-\sqrt{\tau}d_{\partial\Omega}(x)}$ as
$\tau\longrightarrow\infty$.  The estimates are stated in
subsection 2.3 and their proof is given in subsection 2.4. It is
based on the integral representation of $v_g$ with a single layer
potential over $\partial\Omega$. The proof of Theorems 1.2, 1.3
and 1.4 can be done along with the same line as \cite{IK2} in which the
case when $\partial D$ is perfectly insulated is considered.
For reader's convenience we describe an outline of the
proof in section 3. In Appendix we give detailed proofs of four
claims used in subsection 2.4.

\section{Extracting depth}

\subsection{Preliminaries about the direct problem}
In this subsection, following \cite{DL} we describe what we mean by the solution (1.1).
The presentation here is almost parallel to subsection 2.1 in \cite{IK2}.

We put
$\displaystyle
W(0,\,T;H^1(\Omega),(H^1(\Omega))')
=\{u\in L^2(0,\,T;H^1(\Omega))\,\vert\,u'\in L^2(0,\,T;(H^1(\Omega))')\}$.
Given $f\in L^2(0,\,T;H^{-1/2}(\partial\Omega))$ we say that $u\in W(0,\,T;H^1(\Omega),
(H^1(\Omega))')$ satisfy
$$\begin{array}{c}
\displaystyle
\partial_tu-\nabla\cdot\gamma\nabla u=0\,\,\text{in}\,\Omega\times]0,\,T[,\\
\\
\displaystyle
\gamma\nabla u\cdot\nu=f\,\,\text{on}\,\partial\Omega\times]0,\,T[
\end{array}
\tag {2.1}
$$
in the weak sense if $u$ satisfies
$$
\displaystyle
<u'(t),\varphi>
+\int_{\Omega}\gamma(x)\nabla u(x,t)\cdot\nabla\varphi(x)dx
=<f(t),\,\varphi\vert_{\partial\Omega}>\,\,\text{in}\,(0,\,T),
\tag {2.2}
$$
in the sense of distribution on $(0,\,T)$ for all $\varphi\in H^1(\Omega)$ and a.e. $t\in]0,\,T[$.
We see that every $u\in W(0,\,T;H^1(\Omega), (H^1(\Omega))')$
is almost everywhere equal to
a continuous function of $[0,\,T]$ in $L^2(\Omega)$ (Theorem 1 on p.473 in \cite{DL}).  Further, we have
$$\displaystyle
W(0,\,T;H^1(\Omega), (H^1(\Omega))')\hookrightarrow
C^0([0,\,T];L^2(\Omega)),
\tag {2.3}
$$
the space $C^0([0,\,T];L^2(\Omega))$ being equipped with the norm
of uniform convergence. Thus one can consider $u(0)$ and $u(T)$ as
elements of $L^2(\Omega)$.
Then we see that given $u_0\in L^2(\Omega)$ there
exists unique $u$ such that $u$ satisfies (2.1) in the weak
sense and satisfies the initial condition $u(0)=u_0$ (Theorems 1 and 2 on p.512 in \cite{DL}).

Let $u_0=0$.  Remark 2 on p.512 and Theorem 3 on p.520 in
\cite{DL} yields the continuity of $u$ on $f$: there exists
$C_T>0$ independent of $f$ such that
$$\displaystyle
\Vert u\Vert_{L^2(0,T;H^1(\Omega))}
\le C_T\Vert f\Vert_{L^2(0,T;H^{-1/2}(\partial\Omega))}.
\tag {2.4}
$$
Moreover, from (2.2) and (2.4) we have
$$\displaystyle
\Vert u'\Vert_{L^2(0,T;H^1(\Omega)')}
\le C_T\Vert f\Vert_{L^2(0,T;H^{-1/2}(\partial\Omega))}.
$$
This together with (2.3) and (2.4) yields one of the important estimates
in the enclosure method:
$$\displaystyle
\Vert u(T)\Vert_{L^2(\Omega)}
\le C_T\Vert f\Vert_{L^2(0,T;H^{-1/2}(\partial\Omega))}.
\tag {2.5}
$$
In the following subsection, we denote by $u_f$ the weak solution
of (2.1) with $u(0)=0$ and this is the meaning of the weak
solution of (1.1).

\subsection{A basic identity}

Define
$$\displaystyle
w_f(x;\tau)=\int_0^T e^{-\tau t}u_f(x,t)dt,\,\,x\in\Omega
$$
and
$$\displaystyle
g_f(x;\tau)=\int_0^Te^{-\tau t}f(x,t)dt,\,\,x\in\partial\Omega,
$$
where $\tau>0$ is a parameter.  This type of transform has been used
in the study \cite{I4} for the corresponding problem in a one-space
dimensional case.

In this subsection, we derive an identity that connects the data for the parabolic equation
with the Cauchy data of the solutions of the modified Helmholtz-type equations.

$w_f=w$ satisfies
$$\begin{array}{c}
\displaystyle
(\nabla\cdot\gamma\nabla-\tau)w=e^{-\tau T}u_f(x,T)\,\,\text{in}\,\Omega,\\
\\
\displaystyle
\gamma\nabla w\cdot\nu=g_f\,\,\text{on}\,\partial\Omega.
\end{array}
$$

Let $v=v(x)$ satisfy $(\triangle-\tau)v=0$ in $\Omega$.
Integration by parts yields
$$\displaystyle
\int_{\partial\Omega}\left(g_fv-w_f\frac{\partial v}{\partial\nu}\right)dS
=\int_{\Omega}(\gamma-I_3)\nabla v\cdot\nabla w_f dx
+e^{-\tau T}\int_{\Omega}u_f(x,T)v(x)dx.
\tag {2.6}
$$
Let $p_f=p$ be the unique solution of the boundary value problem:
$$\begin{array}{c}
\displaystyle
(\nabla\cdot\gamma\nabla-\tau)p=0\,\,\text{in}\,\Omega,\\
\\
\displaystyle
\gamma\nabla p\cdot\nu=g_f\,\,\text{on}\,\partial\Omega.
\end{array}
$$
Set $\epsilon_f=w_f-p_f$.
Since we have
$$\displaystyle
\int_{\partial\Omega}\left(g_fv-p_f\frac{\partial v}{\partial\nu}\right)dS=
\int_{\Omega}(\gamma-I_3)\nabla v\cdot\nabla p_f dx,
$$
from (2.6) it follows that
$$\begin{array}{c}
\displaystyle
\int_{\partial\Omega}\left(g_fv-w_f\frac{\partial v}{\partial\nu}\right)dS
=\int_{\partial\Omega}\left(g_fv-p_f\frac{\partial v}{\partial\nu}\right)dS
\\
\\
\displaystyle
+\int_{\Omega}(\gamma-I_3)\nabla v\cdot\nabla\epsilon_fdx
+e^{-\tau T}\int_{\Omega}u_f(x,T)v(x)dx.
\end{array}
\tag {2.7}
$$
Note that $\epsilon_f=\epsilon$ satisfies
$$\begin{array}{c}
\displaystyle
(\nabla\cdot\gamma\nabla-\tau)\epsilon=e^{-\tau T}u_f(x,T)\,\,\text{in}\,\Omega,\\
\\
\displaystyle
\gamma\nabla\epsilon\cdot\nu=0\,\,\text{on}\,\partial\Omega.
\end{array}
\tag {2.8}
$$

Let $R_{I_3}(\tau)$ and $R_{\gamma}(\tau)$ denote the Neumann-to-Dirichlet maps
on $\partial\Omega$ for the operators
$\triangle-\tau$ and $\nabla\cdot\gamma\nabla-\tau$ in $\Omega$, respectively.
We have
$$\displaystyle
R_{I_3}(\tau)\left(\frac{\partial v}{\partial\nu}\vert_{\partial\Omega}\right)=v\vert_{\partial\Omega},
\,\,
R_{\gamma}(\tau)g_f
=p_f\vert_{\partial\Omega}.
$$
Since both $R_{I_3}(\tau)$ and $R_{\gamma}(\tau)$ are symmetric,
we obtain from (2.7)
$$\begin{array}{c}
\displaystyle
\int_{\partial\Omega}\left(g_fv-w_f\frac{\partial v}{\partial\nu}\right)dS
=\int_{\partial\Omega}g_f(R_{I_3}(\tau)-R_{\gamma}(\tau))\left(\frac{\partial v}{\partial\nu}\vert_{\partial\Omega}\right)dS
\\
\\
\displaystyle
+\int_{\Omega}(\gamma-I_3)\nabla v\cdot\nabla\epsilon_fdx
+e^{-\tau T}\int_{\Omega}u_f(x,T)v(x)dx.
\end{array}
\tag {2.9}
$$
This is our basic identity.  In the proof of Theorems 1.1 to 1.4
we show that, in some sense, one can ignore the second and third terms of this
right-hand side.  Thus the basic identity provides us
a relationship between the boundary data for the parabolic equation
over a finite time interval and the Cauchy data for the modified Helmholtz-type
equations.

\subsection{Proof of Theorem 1.1}

Since the $\epsilon$ satisfies (2.8), one gets
$$\displaystyle
\Vert\nabla\epsilon\Vert_{L^2(\Omega)}
\le Ce^{-\tau T}\Vert u_f(\,\cdot\,,T)\Vert_{L^2(\Omega)},
\tag {2.10}
$$
where $C$ is a positive constant.  Since $f$ is independent of $\tau$, it follows from (2.5)
and (2.10) that $\Vert\nabla\epsilon\Vert_{L^2(\Omega)}=O(e^{-\tau T})$ as $\tau\longrightarrow\infty$.

Now substitute $v=v_g$ into (2.9). From (1.2) and (1.3) one gets
$\Vert v_g(\,\cdot\,;\tau)\Vert_{H^1(\Omega)}=O(\tau^{-\mu})$ as
$\tau\longrightarrow\infty$. From these one gets the estimate on
the second and third term in (2.9) as $\tau\longrightarrow\infty$:
$$\displaystyle
\int_{\Omega}(\gamma-I_3)\nabla v_g\cdot\nabla\epsilon_fdx
+e^{-\tau T}\int_{\Omega}u_f(x,T)v_gdx
=O(\tau^{-\mu}e^{-\tau T}).
$$
Note that this is a very rough estimate; however, for our purpose
it is enough; at this step we never make use of the assumption
that $\gamma(x)=I_3$ outside $D$.

Summing up, we have obtained the asymptotic formula:
$$\begin{array}{c}
\displaystyle
\int_{\partial\Omega}\left(gv_g-w_fg\right)dS
=\int_{\partial\Omega}g(R_{I_3}(\tau)-R_{\gamma}(\tau))gdS
+O(\tau^{-\mu}e^{-\tau T}).
\end{array}
\tag {2.11}
$$
The following system of inequalities is quite useful to give an estimation of
the first term of the right-hand side.

\proclaim{\noindent Proposition 2.1.}
Let $\gamma_0$ and $\gamma$ satisfy (G1) and (G2).
Let $u$ solve
$$\begin{array}{c}
\displaystyle
\nabla\cdot\gamma\nabla u-\tau u=0\,\,\text{in}\,\Omega,\\
\\
\displaystyle
\gamma\nabla u\cdot\nu=g\,\,\text{on}\,\partial\Omega
\end{array}
$$
and $v$
$$\begin{array}{c}
\displaystyle
\nabla\cdot\gamma_0\nabla v-\tau v=0\,\,\text{in}\,\Omega,\\
\\
\displaystyle
\gamma_0\nabla v\cdot\nu=g\,\,\text{on}\,\partial\Omega.
\end{array}
$$
Then it holds that
$$\displaystyle
\int_{\Omega}(\gamma_0^{-1}-\gamma^{-1})\gamma_0\nabla v\cdot\gamma_0\nabla vdx
\le\int_{\partial\Omega}g(v-u)dS
\le\int_{\Omega}(\gamma-\gamma_0)\nabla v\cdot\nabla v dx.
\tag {2.12}
$$
\endproclaim

For the proof see \cite{IES}.
In the present situation $\gamma_0(x)\equiv I_3$ and $\gamma(x)=I_3$ a.e. $x\in\Omega\setminus D$
and thus from (2.12) we obtain
$$\displaystyle
\int_{D}(I_3-\gamma^{-1})\nabla v_g\cdot\nabla v_gdx
\le
\int_{\partial\Omega}g(R_{I_3}(\tau)-R_{\gamma}(\tau))gdS
\le\int_{D}(\gamma-I_3)\nabla v_g\cdot\nabla v_g dx.
\tag {2.13}
$$

Here we describe a key lemma whose proof is given in the next subsection.
\proclaim{\noindent Lemma 2.1.}
There exist real numbers $\lambda_1$ and $\lambda_2$ independent of $\tau$
such that
$$\displaystyle
\limsup_{\tau\longrightarrow\infty}\tau^{\lambda_1}
e^{2\sqrt{\tau}\text{dist}\,(D,\partial\Omega)}
\int_D\vert\nabla v_g\vert^2dx<\infty
\tag {2.14}
$$
and
$$\displaystyle
\liminf_{\tau\longrightarrow\infty}\tau^{\lambda_2}e^{2\sqrt{\tau}\text{dist}\,(D,\partial\Omega)}\int_D\vert\nabla v_g\vert^2dx>0.
\tag {2.15}
$$

\endproclaim

From the proof one can choose $\lambda_1=2\mu-1$ and $\lambda_2=2\mu+5/2$; however,
the exact values of $\lambda_1$, $\lambda_2$ are not important for the derivation of formula (1.4) itself.

From (2.11), (2.13) and (2.14) one gets
$$\displaystyle
\limsup_{\tau\longrightarrow\infty}
\tau^{\lambda_1}e^{2\sqrt{\tau}\text{dist}\,(D,\partial\Omega)}
\left\vert\int_{\partial\Omega}\left(gv_g-w_fg\right)dS\right\vert<\infty.
\tag {2.16}
$$

Now consider the case when (A1) is satisfied.  It follows from the right half of (2.13) and (2.15) that
$$\displaystyle
\liminf_{\tau\longrightarrow\infty}\tau^{\lambda_2}
e^{2\sqrt{\tau}\text{dist}(D,\partial\Omega)}\left(-\int_{\partial\Omega}g(R_{I_3}(\tau)-R_{\gamma}(\tau))gdS\right)
>0.
$$
This together with (2.11) gives
$$\displaystyle
\liminf_{\tau\longrightarrow\infty}\tau^{\lambda_2}
e^{2\sqrt{\tau}\text{dist}(D,\partial\Omega)}\left(
-\int_{\partial\Omega}\left(gv_g-w_fg\right)dS \right)
>0.
\tag {2.17}
$$
This also implies also that there exists $\tau_0>0$ such that for all $\tau\ge\tau_0$
$$\displaystyle
-\int_{\partial\Omega}\left(gv_g-w_fg\right)dS>0.
$$

Next consider the case when (A2) is satisfied.  Since $I_3-\gamma(x)^{-1}=\gamma(x)^{-1/2}
(\gamma(x)-I_3)\gamma(x)^{-1/2}$, one can find a positive constant $C$ such that,
for all $\xi\in\Bbb R^3$ $\displaystyle (I_3-\gamma(x)^{-1})\xi\cdot\xi\ge C\vert\xi\vert^2$.
Hence a similar argument yields that
$$\displaystyle
\liminf_{\tau\longrightarrow\infty}\tau^{\lambda_2}
e^{2\sqrt{\tau}\text{dist}(D,\partial\Omega)}\left(
\int_{\partial\Omega}\left(gv_g-w_fg\right)dS \right)
>0
\tag {2.18}
$$
and this implies that there exists $\tau_0>0$ such that for all $\tau\ge\tau_0$
$$\displaystyle
\int_{\partial\Omega}\left(gv_g-w_fg\right)dS>0.
$$

Now formula (1.4) is a direct consequence of (2.16), (2.17), (2.18)
and the identity
$$\displaystyle
\int_{\partial\Omega}\int_0^Te^{-\tau t}\left(v(x)f(x,t;\tau)-u_f(x,t)\frac{\partial v}{\partial\nu}(x)\right)dtdS
=\int_{\partial\Omega}\left(g_fv-w_f\frac{\partial v}{\partial\nu}\right)dS.
$$

\subsection{Proof of Lemma 2.1}
Let $\mu$ be the constant in (1.2).  Set $\tilde{v}(x;\tau)=\tau^{\mu}v_g(x;\tau)$
and
$$\displaystyle
\tilde{g}(x;\tau)=\tau^{\mu}\int_0^Te^{-\tau t}f(x,t)dt,\,\,x\in\partial\Omega.
$$
It suffices to prove (2.14) and (2.15) for $\tilde{v}$ instead of $v_g$.

$\tilde{v}$ satisfies $(\triangle-\tau)\tilde{v}=0$ in $\Omega$ and
$\partial\tilde{v}/\partial\nu=\tilde{g}$ on $\partial\Omega$.
In what follows we simply write $\tilde{v}$ and $\tilde{g}$ as $v$ and $g$, respectively.
We think that this makes no confusion.
Thus from (1.2) one has the following: there exist positive constants
$C$ and $\tau_0$ independent of $x\in\partial\Omega$ such that, for all $x\in\partial\Omega$
and all $\tau\ge\tau_0$,
$$\displaystyle
C^{-1}\le g(x;\tau)\le C.
\tag {2.19}
$$

Using the potential theory (cf \cite{Mi}), one has the expression
$$\displaystyle
v(x;\tau)=\frac{1}{2\pi}
\int_{\partial\Omega}\frac{e^{-\sqrt{\tau}\vert x-y\vert}}{\vert x-y\vert}\psi(y;\tau)dS_y,\,\,x\in\Omega,
$$
where $\psi(\,\cdot\,;\tau)\in C(\partial\Omega)$ is the unique solution of the integral equation of
the second kind on $\partial\Omega$:
$$\displaystyle
\psi(y;\tau)+\frac{1}{2\pi}\int_{\partial\Omega}\frac{\partial}{\partial\nu_y}\left(\frac{e^{-\sqrt{\tau}\vert y-y'\vert}}
{\vert y-y'\vert}\right)\psi(y';\tau)dS_{y'}=g(y;\tau),\,\,y\in\partial\Omega.
\tag {2.20}
$$
It is well known that the operator
$$\displaystyle
C(\partial\Omega)\ni\varphi\longmapsto S_{\partial\Omega}(\tau)\varphi\in C(\partial\Omega),
$$
where
$$
\displaystyle
S_{\partial\Omega}(\tau)\varphi(y)
=\frac{1}{2\pi}\int_{\partial\Omega}\frac{\partial}{\partial\nu_y}\left(\frac{e^{-\sqrt{\tau}\vert y-y'\vert}}
{\vert y-y'\vert}\right)\varphi(y')dS_{y'},\,\,y\in\partial\Omega
$$
is bounded and its operator norm has a bound $O(\tau^{-1/2})$ as $\tau\longrightarrow\infty$.
Thus it follows from (2.19)
and (2.20) that $\psi(\,\cdot\,;\tau)$ also has the following: there exist positive constants
$C'$ and $\tau_0$ independent of $x\in\partial\Omega$ such that, for all $x\in\partial\Omega$
and all $\tau\ge\tau_0$
$$\displaystyle
C'^{-1}\le \psi(x;\tau)\le C'.
\tag {2.21}
$$
Since
$$\displaystyle
\nabla v(x;\tau)
=-\frac{1}{2\pi}
\int_{\partial\Omega}\psi(y;\tau)
\frac{e^{-\sqrt{\tau}\vert x-y\vert}}
{\vert x-y\vert}
\left(\sqrt{\tau}+\frac{1}{\vert x-y\vert}\right)\frac{x-y}{\vert x-y\vert}dS_y,
$$
one has
$$\displaystyle
\int_D\vert\nabla v(x;\tau)\vert^2dx
=\int_{\partial\Omega}dS_y\int_{\partial\Omega}dS_{y'}
\int_D dx
e^{-\sqrt{\tau}\vert x-y\vert}
e^{-\sqrt{\tau}\vert x-y'\vert}
\Phi(x,y,y';\tau),
\tag {2.22}
$$
where
$$\begin{array}{c}
\displaystyle
\Phi(x,y,y';\tau)\\
\\
\displaystyle
=\frac{1}{(2\pi)^2}
\frac{\psi(y;\tau)\psi(y';\tau)}
{\vert x-y\vert\vert x-y'\vert}
\left(\sqrt{\tau}+\frac{1}{\vert x-y\vert}\right)
\left(\sqrt{\tau}+\frac{1}{\vert x-y'\vert}\right)
\frac{x-y}{\vert x-y\vert}\cdot\frac{x-y'}{\vert x-y'\vert}.
\end{array}
\tag {2.23}
$$
From (2.21), (2.22) and (2.23) we can easily obtain (2.14) with $\lambda_1=-1$.

The problem is the proof of (2.15).
We divide the integrand of (2.22) into two parts.
Set $d_0=\text{dist}\,(D,\partial\Omega)$ and
${\cal M}=\{(x,y)\in\overline D\times\partial\Omega\,\vert\,\vert x-y\vert=d_0\}$.
It is easy to see that ${\cal M}$ coincides with the set of all $(x,y)\in\partial D\times\partial\Omega$
such that $\vert x-y\vert=d_0$.

In what follows we denote by $B_{R}(z)$ the open ball centered at a point $z$ with radius $R$.
Given $\delta>0$ define
$$\displaystyle
{\cal W}_{\delta}=\cup_{(x_0,y_0)\in{\cal M}}(\overline D\cap B_{\delta}(x_0))
\times(\partial\Omega\cap B_{\delta}(y_0))\times(\partial\Omega\cap B_{\delta}(y_0)).
$$
The set ${\cal W}_{\delta}$ is open in $\overline D\times\partial\Omega\times\partial\Omega$ and
contains the set of all $(x,y,y)$ with $(x,y)\in{\cal M}$.

Here we state two claims concerning the ${\cal W}_{\delta}$ whose
proof is given in Appendix.

\noindent
{\bf Claim 1.}  Given $\epsilon>0$ there exists $\delta_1>0$ such that
for all $(x,y,y')\in {\cal W}_{\delta_1}$ it holds that
$$\displaystyle
\frac{x-y}{\vert x-y\vert}\cdot\frac{x-y'}{\vert x-y'\vert}
\ge 1-\epsilon,\,\,
\vert x-y\vert\le d_0+\epsilon,\,\,
\vert x-y'\vert\le d_0+\epsilon.
$$

\noindent
{\bf Claim 2.}  Given $\delta_1>0$ there exists $\delta_2>0$ such that
if $(x,y,y')\in\overline D\times\partial\Omega\times\partial\Omega\setminus{\cal W}_{\delta_1}$,
then $\vert x-y\vert+\vert x-y'\vert\ge 2d_0+\delta_2$.

Thus giving $\epsilon=1/2$ in claim 1 and choosing the corresponding $\delta_1$,
we have the following: if $(x,y,y')\in {\cal W}_{\delta_1}$, then
$$\displaystyle
\frac{x-y}{\vert x-y\vert}\cdot\frac{x-y'}{\vert x-y'\vert}
\ge\frac{1}{2}.
$$
It follows from this together with the left half of (2.21) and (2.23) that there exist constants $C_1$
and $\tau_0>$ such that, for all $(x,y,y')\in {\cal W}_{\delta_1}$ and $\tau\ge\tau_0$,
$$\displaystyle
\Phi(x,y,y';\tau)\ge C_1\tau.
$$
On the other hand, using the right half of (2.21), it is easy to see that
there exist positive constants $C_2$ and $\tau_1>\tau_0$ such that, for all
$(x,y,y')\in\overline D\times\partial\Omega\times\partial\Omega$ and $\tau\ge\tau_1$,
$$\displaystyle
\Phi(x,y,y';\tau)\le C_2\tau.
$$

Now choose $\delta_2$ in claim 2 corresponding to $\delta_1$ already chosen.

Then we have
$\displaystyle e^{-\sqrt{\tau}(\vert x-y\vert+\vert x-y'\vert)} \le e^{-2\sqrt{\tau}d_0-\sqrt{\tau}\delta_2}$
for any $(x,y,y')\in\overline D\times\partial\Omega\times\partial\Omega\setminus{\cal W}_{\delta_1}$.
Hence dividing the integral (2.22) into ${\cal W}_{\delta_1}$ and its complement, one
gets as $\tau\longrightarrow\infty$
$$\displaystyle
\tau^{-1}\int_{D}\vert\nabla v(x;\tau)\vert^2dx
\ge C_1\int_{{\cal W}_{\delta_1}}dS_ydS_{y'}dx e^{-\sqrt{\tau}(\vert x-y\vert+\vert x-y'\vert)}
+O(e^{-2\sqrt{\tau}d_0-\sqrt{\tau}\delta_2}).
\tag {2.24}
$$
Choose $(x_0,y_0)\in{\cal M}$.
It follows from the definition of ${\cal W}_{\delta_1}$ and the inequality $\vert x-y\vert+\vert x-y'\vert
\le 2\vert x-x_0\vert+2d_0+\vert y_0-y\vert+\vert y_0-y'\vert$ that
$$\begin{array}{c}
\displaystyle
\int_{{\cal W}_{\delta_1}}dS_ydS_{y'}dx e^{-\sqrt{\tau}(\vert x-y\vert+\vert x-y'\vert)}\\
\\
\displaystyle
\ge
\int_{\overline D\cap B_{\delta_1}(x_0)}dx
\int_{\partial\Omega\cap B_{\delta_1}(y_0)}dS_y\int_{\partial\Omega\cap B_{\delta_1}(y_0)}dS_{y'}
e^{-\sqrt{\tau}(\vert x-y\vert+\vert x-y'\vert)}\\
\\
\displaystyle
\ge
e^{-2\sqrt{\tau}d_0}
\int_{\overline D\cap B_{\delta_1}(x_0)}e^{-2\sqrt{\tau}\vert x-x_0\vert}dx
\left(\int_{\partial\Omega\cap B_{\delta_1}(y_0)} e^{-\sqrt{\tau}\vert y_0-y\vert}dS_y\right)^2.
\end{array}
$$

Now (2.15) with $\lambda_2=5/2$ is a direct consequence of this together with (2.24) and the following two claims.

\noindent
{\bf Claim 3.} For all $\delta>0$ we have
$$\displaystyle
\liminf_{\tau\longrightarrow\infty}
(\sqrt{\tau})^3\int_{D\cap B_{\delta}(x_0)}e^{-2\sqrt{\tau}\vert x-x_0\vert}dx>0.
$$

\noindent
{\bf Claim 4.}  For all $\delta>0$ we have
$$\displaystyle
\liminf_{\tau\longrightarrow\infty}
\tau\int_{\partial\Omega\cap B_{\delta}(y_0)}e^{-\sqrt{\tau}\vert y_0-y\vert}dS_y>0.
$$

\noindent
For the proof of these claims see Appendix.

\section{Outline of the proof of Theorems 1.2, 1.3 and 1.4}

Since from (1.5) we have
$$\displaystyle
g_f(x;\tau)=\frac{\partial v}{\partial\nu}(x;\tau)\int_0^Te^{-\tau t}\varphi(t)dt,
$$
it follows from (2.9) that
$$\begin{array}{c}
\displaystyle
\int_{\partial\Omega}\left(g_fv-w_f\frac{\partial v}{\partial\nu}\right)dS
=\int_0^Te^{-\tau t}\varphi(t)dt\int_{\partial\Omega}\frac{\partial v}{\partial\nu}(R_{I_3}(\tau)-R_{\gamma}(\tau))
\left(\frac{\partial v}{\partial\nu}\vert_{\partial\Omega}\right)dS
\\
\\
\displaystyle
+\int_{\Omega}(\gamma-I_3)\nabla v\cdot\nabla\epsilon_fdx
+e^{-\tau T}\int_{\Omega}u_f(x,T)v(x;\tau)dx.
\end{array}
\tag {3.1}
$$
By virtue of (2.5) one knows that for both $v$ in Theorems 1.2-1.4 there exists a constant
$\kappa$ such that $\Vert u_f(\,\cdot\,,;T)\Vert_{L^2(\Omega)}=O(e^{\kappa\sqrt{\tau}})$ as
$\tau\longrightarrow\infty$.  This together (2.10) yields
that as $\tau\longrightarrow\infty$
$$\displaystyle
\int_{\Omega}(\gamma-I_3)\nabla v\cdot\nabla\epsilon_fdx
+e^{-\tau T}\int_{\Omega}u_f(x,T)v(x;\tau)dx
=O(e^{-\tau T/2}).
\tag {3.2}
$$
Here we recall the following lemma.
\proclaim{\noindent Lemma 3.1.}
There exist real numbers $\mu_1, \mu_2, \mu_3\in\Bbb R$
such that, for all $\omega\in S^2$, $p\in\Bbb R^3\setminus\overline\Omega$ and $y\in\Bbb R^3$,
$$\displaystyle
\liminf_{\tau\longrightarrow\infty}
\tau^{\mu_1}e^{-2\sqrt{\tau}h_D(\omega)}\int_D e^{2\sqrt{\tau}x\cdot\omega} dx>0,
\tag {3.3}
$$
$$\displaystyle
\liminf_{\tau\longrightarrow\infty}
\tau^{\mu_2}e^{2\sqrt{\tau}d_D(p)}\int_D e^{-2\sqrt{\tau}\vert x-p\vert} dx>0
\tag {3.4}
$$
and
$$\displaystyle
\liminf_{\tau\longrightarrow\infty}
\tau^{\mu_3}e^{-2\sqrt{\tau}R_D(y)}\int_D e^{2\sqrt{\tau}\vert x-y\vert} dx>0.
\tag {3.5}
$$

\endproclaim

It is not necessary for us to know the values
of $\mu_1$, $\mu_2$, $\mu_3$ precisely.
Every case can be reduced to the case when $D$ is given by an open ball since we are assuming that $\partial D$ is smooth.
See \cite{IE} for the proof of (3.3) and \cite{IK2}(or \cite{I}) for the proof (3.4) and (3.5).
Now it is a due course to see that a combination
of (1.6), (2.12), (3.1), (3.2) and (3.3)/(3.4)/(3.5) yields (1.7)/(1.9)/(1.11).

\section{Conclusion and open problems}

We showed how the enclosure method can be applied to
inverse initial boundary value problems over a finite time interval
for parabolic equations with discontinuous coefficients.
We established four types of formulae.
It should be emphasized that in all formulae the initial temperature
field {\it inside the body} is assumed to be a {\it known
constant}. We think that this is a natural condition and can be
realized without special care in practice. In fact, just make it
cold by using a refrigerator if the size of the body is not so
large!

Two of them are {\it new} in idea and yield the following: (I) a {\it depth} of unknown
inclusions in a heat conductive body from the surface of the body
with a single set of a heat flux and the corresponding temperature
on the surface over a finite time interval (see Theorem 1.1);
(II) the {\it minimum radius} of the open ball centered at a given point that contains unknown inclusions
with a special explicit heat flux with a large parameter (see Theorem 1.4).

The point is the choice of the heat flux $f$ and a solution $v$ of the modified Helmholtz equation
$(\triangle-\tau)v=0$ in $\Omega$
in the integral
$$\displaystyle
\int_{\partial\Omega}\int_0^Te^{-\tau t}\left(vf-u_f\frac{\partial v}{\partial\nu}\right)dtdS.
\tag {4.1}
$$
In (I) first $f$ is given and we choose $v$ by {\it solving} the Neumann problem for the modified Helmholtz equation in $\Omega$
whose Neumann data can be calculated from $f$; in (II) using a special $v$ which is {\it growing everywhere} as $\tau\longrightarrow\infty$, we specify the form of $f$.

The procedure suggested from (I) of extracting $\text{dist}\,(D,\partial\Omega)$ is extremely simple and summarized as follows.

\noindent
(i)  Give the heat flux $f$ satisfying (1.2) for $\mu\in\Bbb R$
across $\partial\Omega$ over the time interval $]0,\,T[$.

\noindent
(ii) Measure the temperature $u_f(x,t)$ on $\partial\Omega$ over the time interval $]0,\,T[$.

\noindent
(iii)  Fix a large $\tau>0$ and compute the solution $v_g$ of (1.3).

\noindent
(iv)  Compute the quantity
$$\displaystyle
-\frac{1}{2\sqrt{\tau}}
\log\left\vert\int_{\partial\Omega}
\int_0^Te^{-\tau t}\left(v(x;\tau)f(x,t)-u_f(x,t)\frac{\partial v}{\partial\nu}(x;\tau)\right)dtdS\right\vert
$$
as an approximation of $\text{dist}\,(D,\partial\Omega)$.

If $D$ is near surface $\partial\Omega$ and isolated in a small part,
the information $\text{dist}\,(D,\partial\Omega)$ may not be so useful; however,
if $D$ is deep inside or occupies a large part, then the set of all $x\in\Omega$ such that
$\text{dist}\,(D,\partial\Omega)<d_{\partial\Omega}(x)$ may give a good estimation of $D$ from above.

The method can be applied also to more
complicated situations, for example, inclusions in a body with a known {\it inhomogeneous} isotropic or anisotropic
conductivity apart from some technical difficulties or similar problems with {\it acoustic/elastic/electromagnetic waves}, etc..
Such applications belong to our future
study.

A next challenging problem is: to clarify what information about $D$ can be extracted from
the asymptotic behaviour of integral (4.1)
as $\tau\longrightarrow\infty$ if $f$ is {\it fixed};
$v$ is one of the {\it three special solutions}
of the modified Helmholtz equation in Theorems 1.2, 1.3 and 1.4, that is one of
$v=e^{\sqrt{\tau}\,x\cdot\omega}$, (1.8) and (1.10).
This remains open at the present time.

Finally it should be pointed out that in Beilina and Klibanov
\cite{BK} an inverse problem for the Cauchy problem for a
hyperbolic equation with a variable smooth coefficient in the
whole space has been considered.  The initial data is given by the
delta function concentrated at a {\it single point} located
outside of the domain of interest. They gave a numerical procedure
for the reconstruction of the coefficient in the domain of
interest by observing the solution on the surface surrounding the
domain.  This procedure is also tested numerically. It would be
interesting to construct and test a numerical method based on the
formulae obtained in our paper.  The study is not in the scope of
our paper, but rather is left for future work.

$$\quad$$

\centerline{{\bf Acknowledgments}}

MI was partially supported by Grant-in-Aid for
Scientific Research (C)(No. 21540162) of Japan  Society for
the Promotion of Science.
MK was partially supported by Grant-in-Aid for
Scientific Research (C)(No. 22540194) of Japan  Society for
the Promotion of Science.
The authors thank the referees and a board member
for many valuable remarks.

$$\quad$$

\section{Appendix. Proof of claims}

\subsection{Claim 1}
Define
$$\displaystyle
F(x,y,y')
=\left(1-\frac{x-y}{\vert x-y\vert}\cdot\frac{x-y'}{\vert x-y'\vert}\right)
+(\vert x-y\vert-d_0)+(\vert x-y'\vert-d_0),\,\,(x,y,y')\in\overline D\times\partial\Omega\times\partial\Omega.
$$
Since $\vert x-y\vert\ge d_0$ and $(x-y)/\vert
x-y\vert\cdot(x-y')/\vert x-y'\vert\le 1$ for
$(x,y,y')\in\overline D\times\partial\Omega\times\partial\Omega$,
it suffices to prove that given $\epsilon>0$ there exists
$\delta_1>0$ such that $F(x,y,y')\le\epsilon$ for all
$(x,y,y')\in{\cal W}_{\delta_1}$. Assume that this is not true.
There exists $\epsilon_0>0$ and a sequence $(x_l,y_l,y'_{l})\in
{\cal W}_{1/l}$, $l=1,2,\cdots$ such that $F(x_l,y_l,y'_{l})\ge
\epsilon_0$.  By the definition of ${\cal W}_{1/l}$ we know that
for each $l$ there exists $(p_l,q_l)\in{\cal M}$ such that
$\vert x_l-p_l\vert<1/l$, $\vert y_l-q_l\vert<1/l$ and $\vert
y'_l-q_l\vert<1/l$. Since $\overline D$ and $\partial\Omega$ are
compact, one can choose a subsequence $l_1,l_2,\cdots$ of
$l=1,2,\cdots$ in such a way that the limits
$\lim_{j\longrightarrow\infty}x_{l_j}=x\in\overline D$,
$\lim_{j\longrightarrow\infty}p_{l_j}=p\in\overline D$,
$\lim_{j\longrightarrow\infty}y_{l_j}=y\in\partial\Omega$,
$\lim_{j\longrightarrow\infty}y'_{l_j}=y'\in\partial\Omega$ and
$\lim_{j\longrightarrow\infty}q_{l_j}=q\in\partial\Omega$ exist.
Clearly it holds that $x=p$ and $y=y'=q$.  Since ${\cal M}$ is
closed, one gets $(x,y)\in{\cal M}$ and thus $\vert x-y\vert=d_0$.
This together with $y=y'$ gives $F(x,y,y')=0$. On the other hand,
since $F(x_{l_j},y_{l_j},y'_{l_j})\ge\epsilon_0$, one gets
$F(x,y,y')\ge\epsilon_0$. This is a contradiction.

\subsection{Claim 2}
Assume that the statement is not true.  There exist a $\delta_0>0$ and
a sequence $(x_l,y_l,y'_l)\in\overline D\times\partial\Omega\times\partial\Omega\setminus {\cal W}_{\delta_0}$,
$l=1,2,\cdots$,
such that $\vert x_l-y_l\vert+\vert x_l-y'_l\vert<2d_0+1/l$.  Since $\overline D$ and $\partial\Omega$ are compact,
if necessary replacing the sequence with a suitable subsequence, one may assume that the limits
$\lim_{l\longrightarrow\infty}x_l=x\in\overline D$, $\lim_{l\longrightarrow\infty}y_l=y\in\partial\Omega$
and $\lim_{l\longrightarrow\infty}y'_l=y'\in\partial\Omega$ exist.  Since ${\cal W}_{\delta_0}$ is open,
one has $(x,y,y')\in\overline D\times\partial\Omega\times\partial\Omega\setminus{\cal W}_{\delta_0}$.

On the other hand, since $(\vert x_l-y_l\vert-d_0) +(\vert
x_l-y'_l\vert-d_0)<1/l$, $\vert x_l-y_l\vert\ge d_0$ and $\vert
x_l-y'_l\vert\ge d_0$, we obtain $\vert x-y\vert=\vert
x-y'\vert=d_0$.  This means that $(x,y)\in{\cal M}$ and
$(x,y')\in{\cal M}$.  Since $x\in\partial D$, using local
coordinates at $x$ and $y$, one can easily show that
$(y-x)/d_0=\nu_x$ and similarly $(y'-x)/d_0=\nu_x$.  This yields
$y=y'$ and thus $(x,y,y')\in{\cal W}_{\delta_0}$. This is a
contradiction.

\noindent
$\Box$

\subsection{Claim 3}
Given $\delta>0$ one can choose $\delta_0>0$ with
$\delta_0<\delta$ such that:there exists a smooth function $g$ on
$\Bbb R^2$ with compact support such that $g(0,0)=0$ and $D\cap
B_{\delta_0}(x_0)=\{x_0+\sigma_1\mbox{\boldmath $e$}_1+
\sigma_2\mbox{\boldmath
$e_2$}-s\nu_{x_0}\,\vert\,\sigma_1^2+\sigma_2^2+s^2<\delta_0^2,
s>g(\sigma_1,\sigma_2)\}$, where $\mbox{\boldmath $e$}_1$ and
$\mbox{\boldmath $e_2$}$ are unit tangent vectors at $x_0$ and
orthogonal each other.

One can choose a positive constant $C$ in such a way that for all $\sigma=(\sigma_1,\sigma_2)\in\Bbb R^2$
it holds that $\vert g(\sigma)\vert\le C\vert\sigma\vert$.
Let $\tau\ge 1$.  We can easily see that if $s>C\vert\sigma\vert$,
then $s>\sqrt{\tau}g(\sigma/\sqrt{\tau})$.  This together with change of variables yields
$$\begin{array}{c}
\displaystyle
(\sqrt{\tau})^3\int_{D\cap B_{\delta_0}(x_0)}e^{-2\sqrt{\tau}\vert x-x_0\vert}dx
\ge
\int_{s>\sqrt{\tau}g(\sigma/\sqrt{\tau}),\,\vert\sigma\vert^2+s^2<(\sqrt{\tau}\delta_0)^2}
e^{-2\sqrt{\vert\sigma\vert^2+s^2}}d\sigma ds\\
\\
\displaystyle
\ge
\int_{s>C\vert\sigma\vert,\,\vert\sigma\vert^2+s^2<(\sqrt{\tau}\delta_0)^2}
e^{-2\sqrt{\vert\sigma\vert^2+s^2}}d\sigma ds.
\end{array}
$$
Since
$$
\displaystyle
\lim_{\tau\longrightarrow\infty}
\int_{s>C\vert\sigma\vert,\,\vert\sigma\vert^2+s^2<(\sqrt{\tau}\delta_0)^2}
e^{-2\sqrt{\vert\sigma\vert^2+s^2}}d\sigma ds
=\int_{s>C\vert\sigma\vert}
e^{-2\sqrt{\vert\sigma\vert^2+s^2}}d\sigma ds<\infty,
$$
we have the desired conclusion.

\subsection{Claim 4}
Given $\delta>0$ one can choose $\delta_0>0$ with
$\delta_0<\delta$ such that there exists a smooth function $h$ on
$\Bbb R^2$ with compact support such that $h(0,0)=0$ and
$\partial\Omega\cap B_{\delta_0}(y_0)=\{y_0+
\sigma_1\mbox{\boldmath $e$}_1+ \sigma_2\mbox{\boldmath
$e$}_2-h(\sigma_1,\sigma_2)\nu_{y_0}\,\vert\,\sigma_1^2+\sigma_2^2+h(\sigma_1,\sigma_2)^2<\delta_0^2\}$,
where $\mbox{\boldmath $e$}_1$ and $\mbox{\boldmath $e_2$}$ are
unit tangent vectors at $y_0$ and orthogonal each other. Note that
$h$ also satisfies that for all
$\sigma=(\sigma_1,\sigma_2)\in\Bbb R^2$ $\vert h(\sigma)\vert\le
C\vert\sigma\vert$, where $C$ is a positive constant. We see that
if $\vert\sigma\vert<\delta_0/(\sqrt{1+C^2})$, then
$\sigma_1^2+\sigma_2^2+h(\sigma_1,\sigma_2)^2<\delta_0^2$.  This
together with a change of variables yields
$$\begin{array}{c}
\displaystyle
\int_{\partial\Omega\cap B_{\delta_0}(y_0)}e^{-\sqrt{\tau}\vert y_0-y\vert}dS_y
\ge
\int_{\vert\sigma\vert<\delta_0/\sqrt{1+C^2}}e^{-\sqrt{\tau}\sqrt{\vert\sigma\vert^2+h(\sigma)^2}}
\sqrt{1+\vert\nabla h(\sigma)\vert^2}d\sigma\\
\\
\displaystyle
\ge\int_{\vert\sigma\vert<\delta_0/\sqrt{1+C^2}}
e^{-\sqrt{\tau}\sqrt{1+C^2}\vert\sigma\vert}d\sigma\\
\\
\displaystyle
=(\sqrt{\tau}\sqrt{1+C^2})^{-2}
\int_{\vert\sigma\vert<\sqrt{\tau}\delta_0}e^{-\vert\sigma\vert}d\sigma.
\end{array}
$$
Now one gets the desired conclusion.

\vskip1cm
\noindent
e-mail address

ikehata@math.sci.gunma-u.ac.jp

kawasita@math.sci.hiroshima-u.ac.jp
\end{document}